\documentclass[11pt]{amsart}
\usepackage{amsthm,amsmath}
\usepackage[sort&compress,numbers]{natbib}
\newtheorem{theorem}{Theorem}
\newtheorem{corollary}{Corollary}
\newtheorem{lemma}{Lemma}

\renewcommand{\baselinestretch}{1.2}
\begin{document}

\newcommand{\ceil}[1]{\ensuremath{\protect\lceil#1\rceil}}
\newcommand{\Oh}[1]{\ensuremath{\protect\mathcal{O}(#1)}}
\newcommand{\bracket}[1]{\ensuremath{\protect\left(#1\right)}}
\newcommand{\quarter}{\ensuremath{\protect\tfrac{1}{4}}}
\newcommand{\half}{\ensuremath{\protect\tfrac{1}{2}}}
\newcommand{\bfrac}[2]{\ensuremath{\protect\bracket{\frac{#1}{#2}}}}

\newcommand{\G}{\ensuremath{\mathcal{G}}}
\newcommand{\PP}{\ensuremath{\mathcal{P}}}
\newcommand{\QQ}{\ensuremath{\mathcal{Q}}}

\newcommand{\SL}{\ensuremath{\mathsf{sl}}}

\newcommand{\e}{\ensuremath{\mathbf{e}}}
\renewcommand{\e}{\textsf{\textup{e}}}

\newcommand{\N}{\ensuremath{\mathbb{N}}}
\newcommand{\Z}{\ensuremath{\mathbb{Z}}}
\newcommand{\R}{\ensuremath{\mathbb{R}}}

\date{\today}
\title[Graphs with Large Geometric Thickness]{Bounded-Degree Graphs have\\ Arbitrarily Large Geometric Thickness}
\author[J.~Bar{\'a}t]{J{\'a}nos Bar{\'a}t}
\address{Bolyai Institute, University of Szeged, Szeged, Hungary}
\email{barat@math.u-szeged.hu}
\author[J.~Matou{\v{s}}ek]{Ji{\v{r}}{\'i} Matou{\v{s}}ek}
\address{Department of Applied Mathematics and Institute for Theoretical Computer Science, Charles University, Prague, Czech Republic}
\email{matousek@kam.mff.cuni.cz}
\author[D.~R.~Wood]{David R. Wood}
\address{Departament de Matem{\'a}tica Aplicada II, Universitat Polit{\`e}cnica de Catalunya, Barcelona, Spain}
\email{david.wood@upc.edu}
\thanks{Research of J.~Bar{\'a}t was supported by a Marie Curie
Fellowship of the European Community under contract number HPMF-CT-2002-01868
and by the OTKA Grant T.49398. Research of D.~Wood is supported by the Government of Spain grant MEC SB2003-0270, and by the projects MCYT-FEDER BFM2003-00368 and Gen.\ Cat 2001SGR00224.}




\begin{abstract} The \emph{geometric thickness} of a graph $G$ is the 
minimum integer $k$ such that there is a straight line drawing of $G$ with 
its edge set partitioned into $k$ plane subgraphs. Eppstein [Separating 
thickness from geometric thickness. In \emph{Towards a Theory of Geometric 
Graphs}, vol.\ 342 of \emph{Contemp.\ Math.}, AMS, 2004] asked whether 
every graph of bounded maximum degree has bounded geometric thickness. We 
answer this question in the negative, by proving that there exists 
$\Delta$-regular graphs with arbitrarily large geometric thickness. In 
particular, for all $\Delta\geq9$ and for all large $n$, there exists a 
$\Delta$-regular graph with geometric thickness at least 
$c\sqrt{\Delta}\,n^{1/2-4/\Delta-\epsilon}$. Analogous results concerning 
graph drawings with few edge slopes are also presented, thus solving open 
problems by Dujmovi{\'c} et~al.\ [Really straight graph drawings. In 
\emph{Proc.\ 12th International Symp.\ on Graph Drawing} (GD '04), vol.\ 
3383 of \emph{Lecture Notes in Comput.\ Sci.}, Springer, 2004] and Ambrus 
et~al.\ [The slope parameter of graphs. Tech.\ Rep.\ MAT-2005-07, 
Department of Mathematics, Technical University of Denmark, 2005]. \end{abstract}


\maketitle

\section{Introduction}
\label{sec:Intro}

A \emph{drawing} of an (undirected, finite, simple) graph represents each vertex by a distinct point in the plane, and represents each edge by a simple closed curve between its endpoints, such that the only vertices an edge intersects are its own endpoints. Two edges \emph{cross} if they intersect at a point other than a common endpoint. A drawing is \emph{plane} if no two edges cross. 

The \emph{thickness} of an (abstract) graph $G$ is the minimum number of planar subgraphs of $G$ whose union is $G$. Thickness is a classical and widely studied graph parameter; see the survey \citep{MOS98}. The \emph{thickness} of a graph drawing $D$ is the minimum number of plane subgraphs of $D$ whose union is $D$. Any planar graph can be drawn with its vertices at prespecified locations \citep{Halton91, PW-GC01}. It follows that a graph with thickness $k$ has a drawing with thickness $k$ \citep{Halton91}. However, in such a representation the edges might be highly curved\footnote{In fact, a polyline drawing of a random perfect matching on $n$ vertices in convex position almost certainly has $\Omega(n)$ bends on some edge \citep{PW-GC01}.}. 

This motivates the notion of geometric thickness, which is a central topic of this paper. A drawing is \emph{geometric}, also called a \emph{geometric graph}, if every edge is represented by a straight line segment. The \emph{geometric thickness} of a graph $G$ is the minimum thickness of a geometric drawing of $G$; see \citep{DEH-JGAA00, Eppstein-AMS, HSV-CGTA99, DEK-SoCG04, DujWoo-GD05}. Geometric thickness was introduced by \citet{Kainen73} under the name \emph{real linear thickness}. 


Consider the relationship between the various thickness parameters and maximum degree. A graph with maximum degree at most $\Delta$ is called \emph{degree-$\Delta$}. \citet{Halton91} proved that the thickness of a degree-$\Delta$ graph is at most $\ceil{\frac{\Delta}{2}}$, and \citet{SSV-IS04} proved that this bound is tight. \citet{DEK-SoCG04} proved that the geometric thickness of a degree-$4$ graph is at most two. \citet{Eppstein-AMS} asked whether graphs of bounded degree have bounded geometric thickness. The first contribution of this paper is to answer this question in the negative. 

\begin{theorem}
\label{thm:GeomThickness}
For all $\Delta\geq9$ and $\epsilon>0$, for all large $n>n(\epsilon)$ and $n\geq c\Delta$, there exists a $\Delta$-regular $n$-vertex graph with geometric thickness at least $$c\sqrt{\Delta}\,n^{1/2-4/\Delta-\epsilon},$$
for some absolute constant $c$.
\end{theorem}

A number of notes on Theorem~\ref{thm:GeomThickness} are in order. \citet{Eppstein-AMS} proved that geometric thickness is not bounded by thickness. In particular, there exists a graph with thickness three and arbitrarily large geometric thickness. Theorem~\ref{thm:GeomThickness} and the above result of \citet{Halton91} imply a similar result (with a shorter proof). Namely, there exists a $9$-regular graph with thickness at most five and with arbitrarily large geometric thickness.

A \emph{book embedding} is a geometric drawing with the vertices in convex position. The \emph{book thickness} of a graph $G$ is the minimum thickness of a book embedding of $G$. Book thickness is also called \emph{page-number} and \emph{stack-number}; see \citep{DujWoo-DMTCS04} for over fifty references on this topic. By definition, the geometric thickness of $G$ is at most the book thickness of $G$. 
On the other hand \citet{Eppstein01} proved that there exists a graph with geometric thickness two and arbitrarily large book thickness
(also see \citep{Blankenship-PhD03, BO99}). Thus book thickness is not bounded by any function of geometric thickness.

Theorem~\ref{thm:GeomThickness} is analogous to a result of \citet{Malitz94a}, who proved that there exists $\Delta$-regular $n$-vertex graphs with book thickness at least $c\sqrt{\Delta} n^{1/2-1/\Delta}$. Malitz's proof is based on a probabilistic construction 
of a graph with certain expansion properties. The proof of Theorem~\ref{thm:GeomThickness} is easily adapted to prove Malitz's result for $\Delta\geq3$. The difference 
in the bounds ($n^{1/2-4/\Delta}$ and $n^{1/2-1/\Delta}$) is caused by the 
difference between the number of order types of point sets in general and 
convex position ($\approx n^{4n}$ and $n!$). \citet{Malitz94a} 
also proved an upper bound of $\Oh{\sqrt{m}}\subseteq\Oh{\sqrt{\Delta n}}$ 
on the book thickness, and thus the geometric thickness, of $m$-edge 
graphs.

The other contributions of this paper concern geometric graph drawings with few slopes. \citet{WadeChu-CJ94} defined the \emph{slope-number} of a graph $G$ to be the minimum number of distinct edge slopes in a geometric drawing of $G$. For example, \citet{Jamison-DM86} proved that the slope-number of $K_n$ is $n$. \citet{ReallyStraight-GD04} asked whether every graph with bounded maximum degree has bounded slope-number. Since edges with the same slope do not cross, geometric thickness is a lower bound on the slope-number. Thus Theorem~\ref{thm:GeomThickness} immediately implies a negative answer to this question for $\Delta\geq9$, which is improved as follows.

\begin{theorem}
\label{thm:SlopeNumber}
For all $\Delta\geq5$ there exists a $\Delta$-regular graph with arbitrarily large slope-number.
\end{theorem}

The proofs of Theorems~\ref{thm:GeomThickness} and \ref{thm:SlopeNumber} are simple. We basically show that there are more graphs with bounded degree than with with bounded geometric thickness (or slope-number). Our counting arguments are based on two tools from the literature that are introduced in Section~\ref{sec:Tools}. Theorem~\ref{thm:GeomThickness} is then proved in Section~\ref{sec:GeometricThickness}. In Section~\ref{sec:SlopeParameter} we study a graph parameter introduced by \citet{ABH-combined} that is similar to the slope-number, and solve two of their open problems. The proofs will also serve as a useful introduction to the proof of Theorem~\ref{thm:SlopeNumber}, which is presented in Section~\ref{sec:SlopeNumber}.

\section{Tools}
\label{sec:Tools}

All of our results are based on the Milnor-Thom theorem \citep{Milnor64, Thom65}, first proved by \citet{PO49}. Let $\PP=(P_1,P_2,\dots,P_t)$ be a system of $d$-variate real polynomials. A vector $\sigma\in\{-1,0,+1\}^n$ is a \emph{sign pattern} of \PP\ if there exists an $x\in\R^d$ such that the sign of $P_i(x)$ is $\sigma_i$, for all $i=1,2,\dots,t$. The following precise version of the Milnor-Thom theorem is due to \citet{PR93}. 

\begin{lemma}[\citep{PR93}]
\label{lem:MilnorThom}
Let $\PP=(P_1,P_2,\dots,P_t)$ be a system of $d$-variate real polynomials of degree at most $D$. Then the number of sign patterns of \PP\ is at most 
$$\bfrac{50\,Dt}{d}^d.$$
\end{lemma}

Some of our proofs only need sign patterns that distinguish between zero and nonzero values. In this setting, \citet{RBG-JAMS01} gave a better bound with a short proof; also see \citep{Mat02}.

Our second tool is a corollary of more precise bounds due to \citet{BC-JCTA78}, \citet{Wormald78}, and \citet{McKay-AC85}; see Appendix~A.

\begin{lemma}[\citep{BC-JCTA78,Wormald78,McKay-AC85}]
\label{lem:NumberRegular}
For all integers $\Delta\geq1$ and $n\geq c\Delta$, 
the number of labelled $\Delta$-regular $n$-vertex graphs is at least 
$$\bfrac{n}{3\Delta}^{\Delta n/2},$$
for some absolute constant $c$.
\end{lemma}

\section{Geometric Thickness: Proof of Theorem~\ref{thm:GeomThickness}}
\label{sec:GeometricThickness}


Observe that a geometric drawing with thickness $k$ can be perturbed so that the vertices are in \emph{general position} (that is, no three vertices are collinear). Thus in this section we can consider point sets in general position without loss of generality. (We cannot make this assumption for drawings with few slopes.)\ 

\begin{lemma}
\label{lem:NumberGeomThickness}
The number of labelled $n$-vertex graphs with geometric thickness at most $k$
is at most $472^{kn}n^{4n+o(n)}$.
\end{lemma}

\begin{proof}
Let $P$ be a fixed set of $n$ labelled points in general position in the plane. \citet{Ajtai82} proved that there are at most $c^n$ plane geometric graphs with vertex set $P$, where $c\leq10^{13}$. \citet{SS-JCTA03} proved that we can take $c=472$. A geometric graph with vertex set $P$ and thickness at most $k$ consists of a $k$-tuple of plane geometric graphs with vertex set $P$. Thus $P$ admits at most $472^{kn}$ geometric graphs with thickness at most $k$. 

Let $P=(p_1,p_2,\dots,p_n)$ and $Q=(q_1,q_2,\dots,q_n)$ be two sets of $n$ points in general position in the plane. Then $P$ and $Q$ have the same \emph{order type} if for all indices $i<j<k$ we turn in the same direction (left or right) when going from $p_i$ to $p_k$ via $p_j$ and when going from $q_i$ to $q_k$ via $q_j$. Say $P$ and $Q$ have the same order type. Then for all $i,j,k,\ell$, the segments $p_ip_j$ and $p_kp_\ell$ cross if and only if $q_iq_j$ and $q_kq_\ell$ cross. Thus $P$ and $Q$ admit the same set of (at most $472^{kn}$) labelled geometric graphs with thickness at most $k$ (when considering $p_i$ and $q_i$ to be labelled $i$). \citet{Alon-Math86} proved (using the Milnor-Thom theorem) that there are at most $n^{4n+o(n)}$ sets of $n$ points with distinct order types. The result follows.
\end{proof}


It is easily seen that Lemmas~\ref{lem:NumberRegular} and \ref{lem:NumberGeomThickness} imply a lower bound of $c(\Delta-8)\log n$ on the geometric thickness of some $\Delta$-regular $n$-vertex graph. To improve this logarithmic bound to polynomial, we now give a more precise analysis of the number of graphs with  bounded geometric thickness that also accounts for the number of edges in the graph.




\begin{lemma}
\label{lem:NumberGeomGraphsEdges}
Let $P$ be a set of $n$ labelled points in general position in the plane. Let 
$g(P,m)$ be the number of $m$-edge plane geometric graphs with vertex set $P$. Then
$$g(P,m)\leq
\begin{cases}
\binom{n}{2m}\cdot 472^{2m}	&\textup{, if } m\leq\frac{n}{2}\\
472^n						&\textup{, if } m>\frac{n}{2}.
\end{cases}
$$
\end{lemma}

\begin{proof}
As in Lemma~\ref{lem:NumberGeomThickness}, $g(P,m)\leq472^n$ regardless of $m$. Suppose that $m\leq\frac{n}{2}$. An $m$-edge graph has at most $2m$ vertices of nonzero degree. Thus every $m$-edge plane geometric graph with vertex set $P$ is obtained by first choosing a $2m$-element subset $P'\subseteq P$, and then choosing a plane geometric graph on $P'$. The result follows.
\end{proof}

\begin{lemma}
\label{lem:NumberGeomThicknessEdges}
Let $P$ be a set of $n$ labelled points in general position in the plane. 
For every integer $t$ such that $\frac{2m}{n}\leq t\leq m$, let $g(P,m,t)$ be the number of $m$-edge geometric graphs with vertex set $P$ and thickness at most $t$. Then $$g(P,m,t)\leq\bfrac{ctn}{m}^{2m},$$ for some absolute constant $c$.
\end{lemma}

\begin{proof}
Fix nonnegative integers $m_1\leq m_2\leq\dots\leq m_t$ such that $\sum_im_i=m$. Let $g(P;m_1,m_2,\dots,m_t)$ be the number of geometric graphs with vertex set $P$ and thickness $t$, such that there are $m_i$ edges in the $i$-th subgraph. Then 
\begin{align*}
g(P;m_1,m_2,\dots,m_t)\leq\prod_{i=1}^tg(P,m_i).
\end{align*}
Now $m_1\leq\frac{n}{2}$, as otherwise $m>\frac{tn}{2}\geq m$. Let $j$ be the maximum index such that $m_j\leq\frac{n}{2}$. By Lemma~\ref{lem:NumberGeomGraphsEdges},
\begin{align*}
g(P;m_1,m_2,\dots,m_t)
&\leq
\bracket{\prod_{i=1}^j\binom{n}{2m_i}472^{2m_i}}
\bracket{472^n}^{t-j}.
\end{align*}
Since $\sum_im_i\leq m-\half(t-j)n$,
\begin{align*}
g(P;m_1,m_2,\dots,m_t)
&\leq
\bracket{\prod_{i=1}^j\binom{n}{2m_i}}
\bracket{472^{2m-(t-j)n}}
\bracket{472^{(t-j)n}}\\
&=
472^{2m} \prod_{i=1}^t\binom{n}{2m_i}.
\end{align*}
We can suppose that $t$ divides $2m$. It follows (see Appendix~B) that 
\begin{align*}
g(P;m_1,m_2,\dots,m_t)
\leq
472^{2m}\binom{n}{2m/t}^t.
\end{align*}
It is well known \citep[Proposition~1.3]{Jukna01} that $\binom{n}{k}<\bfrac{\e n}{k}^k$, where \e\ is the base of the natural logarithm. Thus with $k=2m/t$ we have
\begin{align*}
g(P;m_1,m_2,\dots,m_t)
&<
\bfrac{236\e tn}{m}^{2m}.
\end{align*}
Clearly
\begin{align*}
g(P,m,t)\leq\sum_{m_1,\dots,m_t}g(P;m_1,m_2,\dots,m_t),
\end{align*}
where the sum is taken over all nonnegative integers $m_1\leq m_2\leq\dots\leq m_t$ such that $\sum_im_i=m$. The number of such partitions \citep[Proposition~1.4]{Jukna01} is at most 
$$\binom{t+m-1}{m}<\binom{2m}{m}<2^{2m}.$$
Hence
\begin{align*}
g(P,m,t)\leq 2^{2m}\bfrac{236\e tn}{m}^{2m}
\leq \bfrac{ctn}{m}^{2m}.
\end{align*}
\end{proof}

As in Lemma~\ref{lem:NumberGeomThickness}, we have the following corollary of 
Lemma~\ref{lem:NumberGeomThicknessEdges}.

\begin{corollary}
\label{cor:NumberGeomThicknessEdges}
For all integers $t$ such that $\frac{2m}{n}\leq t\leq m$, the number of labelled $n$-vertex $m$-edge graphs with geometric thickness at most $t$ is at most $$n^{4n+o(n)}\bfrac{ctn}{m}^{2m},$$ for some absolute constant $c$.\qed
\end{corollary}

\begin{proof}[Proof of Theorem~\ref{thm:GeomThickness}] Let $t$ be the minimum integer such that every $\Delta$-regular $n$-vertex graph has geometric thickness at most $t$. Thus the number of $\Delta$-regular $n$-vertex graphs is at most the number of labelled graphs with $\half\Delta n$ edges and geometric thickness at most $t$. By Lemma~\ref{lem:NumberRegular} and Corollary~\ref{cor:NumberGeomThicknessEdges},
$$\bfrac{n}{3\Delta}^{\Delta n/2}
\leq 
n^{4n+o(n)}\bfrac{ct}{\Delta}^{\Delta n}
\leq 
n^{4n+\epsilon n}\bfrac{ct}{\Delta}^{\Delta n},$$
for large $n>n(\epsilon)$. Hence $t\geq\sqrt{\Delta}\,n^{1/2-4/\Delta-\epsilon}/(c\sqrt{3})$.
\end{proof}

It remains open whether geometric thickness is bounded by a constant for graphs with $\Delta\leq8$. The above method is easily modified to prove Malitz's lower bound on book thickness that is discussed in Section~\ref{sec:Intro}.

\begin{theorem}[\citep{Malitz94a}]
\label{thm:BookThickness}
For all $\Delta\geq3$ and $n\geq c\Delta$, there exists a $\Delta$-regular $n$-vertex graph with book thickness at least $$c\sqrt{\Delta}\,n^{1/2-1/\Delta},$$
for some absolute constant $c$.
\end{theorem}

\begin{proof}
Obviously the number of order types for point sets in convex position is $n!$.
As in the proof of Theorem~\ref{thm:GeomThickness}, 
$$\bfrac{n}{3\Delta}^{\Delta n/2}
\leq 
n!\bfrac{ct}{\Delta}^{\Delta n}
\leq 
n^n\bfrac{ct}{\Delta}^{\Delta n}.$$
Hence $t\geq\sqrt{\Delta}\,n^{1/2-1/\Delta}/(c\sqrt{3})$. 
(The constant $c$ can be considerably improved here; for 
example, we can replace $472$ by $16$.)
\end{proof}

\section{The Slope Parameter of \citet{ABH-combined}}
\label{sec:SlopeParameter}

\citet{ABH-combined} introduced the following slope parameter of graphs. Let $P\subseteq\R^2$ be a set of points in the plane. Let $S\subset\R\cup\{\infty\}$ be a set of slopes. Let $G(P,S)$ be the graph with vertex set $P$ where two points $v,w\in P$ are adjacent if and only if the slope of the line $\overline{vw}$ is in $S$. The \emph{slope parameter} of a graph $G$, denoted by $\SL(G)$, is the minimum integer $k$ such that $G\cong G(P,S)$ for some point set $P$ and slope set $S$ with $|S|=k$. Note that $\SL$ is well-defined, since $\SL(G)\leq|E(G)|$. In this section we consider the following two questions of \citet{ABH-combined}:
\begin{itemize}
\item what is the maximum slope parameter of an $n$-vertex graph?
\item do graphs of bounded maximum degree have bounded slope parameter?
\end{itemize}


\begin{lemma}
\label{lem:Basic}
The number of labelled $n$-vertex graphs $G$ with slope parameter $\SL(G)\leq k$ is at most
$$\bfrac{50n^2k}{2n+k}^{2n+k}.$$
\end{lemma}

\begin{proof}  Let $\G_{n,k}$ denote the family of labelled $n$-vertex graphs $G$ with slope parameter $\SL(G)\leq k$. Consider $V(G)=\{1,2,\dots,n\}$ for every $G\in\G_{n,k}$. For every $G\in\G_{n,k}$, there is a point set $P=\{(x_i(G),y_i(G)):1\leq i\leq n\}$ and slope set $S=\{s_\ell(G):1\leq\ell\leq k\}$, such that $G\cong G(P,S)$, where vertex $i$ is represented by the point $(x_i(G),y_i(G))$. Fix one such representation of $G$. Without loss of generality, $x_i(G)\ne x_j(G)$ for distinct $i$ and $j$. Thus every $s_\ell(G)<\infty$. For all $i,j,\ell$ with $1\leq i<j\leq n$ and $1\leq\ell\leq k$, and for every graph $G\in G_{n,k}$, we define the number
$$P_{i,j,\ell}(G):=(y_j(G)-y_i(G))-s_\ell(G)\cdot(x_j(G)-x_i(G)).$$
Consider $$\PP:=\{P_{i,j,\ell}:1\leq i<j\leq n,1\leq\ell\leq k\}$$ to be a 
set of $\binom{n}{2}k$ degree-$2$ polynomials on the set of variables $$\{x_1,x_2,\dots,x_n,y_1,y_2,\dots,y_n,s_1,s_2,\dots,s_k\}.$$ Observe that $P_{i,j,\ell}(G)=0$ if and only if $ij$ is an edge of $G$ and $ij$ has slope $s_\ell$ in the representation of $G$.

Consider two distinct graphs $G,H\in\G_{n,k}$. Without loss of generality, there is an edge $ij$ of $G$ that is not an edge of $H$. Thus $(y_j(G)-y_i(G))-s_\ell(G)\cdot(x_j(G)-x_i(G))=0$ for some $\ell$, and $(y_j(H)-y_i(H))-s_\ell(H)\cdot(x_j(H)-x_i(H))\ne 0$ for all $\ell$. Hence $P_{i,j,\ell}(G)=0\ne P_{i,j,\ell}(H)$. That is, any two distinct graphs in $\G_{n,k}$ are distinguished by the sign of some polynomial in \PP. Hence $|\G_{n,k}|$ is at most the number of sign patterns determined by \PP. By 
Lemma~\ref{lem:MilnorThom} with $D=2$, $t=2n+k$, and $d=\binom{n}{2}k$ we have 
$$|\G_{n,k}|
\leq\bfrac{50\cdot 2\cdot\binom{n}{2}k}{2n+k}^{2n+k}
<\bfrac{50n^2k}{2n+k}^{2n+k}.$$
\end{proof}


In response to the first question of \citet{ABH-combined}, we now prove that there exist graphs with surprisingly large slope parameter. In this paper all logarithms are binary unless stated otherwise.

\begin{theorem}
\label{thm:BigSlope}
For all $\epsilon>0$ and for all sufficiently large $n>n(\epsilon)$, there exists an $n$-vertex graph $G$ with slope parameter
$$\SL(G)\geq\frac{n^2}{(4+\epsilon)\log n}.$$
\end{theorem}

\begin{proof}
Suppose that every $n$-vertex graph $G$ has slope parameter $\SL(G)\leq k$. There are $2^{\binom{n}{2}}$ labelled $n$-vertex graphs. By Lemma~\ref{lem:Basic},
$$2^{\binom{n}{2}}\leq\bfrac{50n^2k}{2n+k}^{2n+k}.$$
For large $n>n(\epsilon)$, 
$$\bracket{50n^2}^{2n+n^2/(4+\epsilon)\log n}<2^{\binom{n}{2}}
\leq
\bfrac{50n^2k}{2n+k}^{2n+k}<\bracket{50n^2}^{2n+k}.$$
Hence $$k>\frac{n^2}{(4+\epsilon)\log n}.$$ The result follows.
\end{proof}


Now we prove that the slope parameter of degree-$\Delta$ graphs is unbounded for $\Delta\geq5$, thus answering the second question of \citet{ABH-combined} in the negative. It remains open whether $\SL(G)$ is bounded for degree-$3$ or degree-$4$ graphs $G$.

\begin{theorem}
For all $\Delta\in\{5,6,7,8\}$, for all $\epsilon$ with $0<\epsilon<\Delta-4$, and for all sufficiently large $n>n(\Delta,\epsilon)$, there exists a $\Delta$-regular $n$-vertex graph $G$ with $$\SL(G)>n^{(\Delta-4-\epsilon)/4}.$$
\end{theorem}

\begin{proof}
Let $k:=n^{(\Delta-4-\epsilon)/4}$. Suppose that for some $\Delta\in\{5,6,7,8\}$, every $\Delta$-regular $n$-vertex graph $G$ has $\SL(G)\leq k$. By Lemmas~\ref{lem:NumberRegular} and \ref{lem:Basic},
$$\bfrac{n}{3\Delta}^{\Delta n/2}
\leq
\bfrac{50n^2k}{2n+k}^{2n+k}
<
\bracket{25nk}^{2n+k}
<
\bracket{25n}^{(\Delta-\epsilon)(2n+k)/4}.$$
For large $n>n(\epsilon,\Delta)$,
$$
\bracket{25n}^{(2\Delta-\epsilon)n/4}
<
\bfrac{n}{3\Delta}^{\Delta n/2}
<
\bracket{25n}^{(\Delta-\epsilon)(2n+k)/4}.$$
Thus 
$2\Delta n-\epsilon n<2\Delta n +\Delta k-2\epsilon n-\epsilon k$. That is, $(\Delta-\epsilon)n^{(\Delta-8-\epsilon)/4}>\epsilon$. Thus $\Delta-8-\epsilon\geq0$ for large $n>n(\Delta,\epsilon)$, which is the desired contradiction for $\Delta\leq8$.
\end{proof}


For $\Delta\geq9$ there are graphs with linear slope parameter.

\begin{theorem}
\label{thm:SlopeParameter}
For all $\Delta\geq9$ and $\epsilon>0$, and for all sufficiently large $n>n(\Delta,\epsilon)$, there exists a $\Delta$-regular $n$-vertex graph $G$ with slope parameter $$\SL(G)>\quarter((1-\epsilon)\Delta-8)n.$$
\end{theorem}

\begin{proof}
Suppose that every $\Delta$-regular $n$-vertex graph $G$ has $\SL(G)\leq\alpha n$ for some $\alpha>0$. By Lemmas~\ref{lem:NumberRegular} and \ref{lem:Basic}, 
$$\bfrac{n}{3\Delta}^{\Delta n/2}
\leq
\bfrac{50\alpha\,n^2}{2+\alpha}^{(2+\alpha)n}.$$
For large $n>n(\Delta,\epsilon)$,
$$\bracket{8\,n}^{(1-\epsilon)\Delta n/2}
<
\bfrac{n}{3\Delta}^{\Delta n/2}
\leq
\bfrac{50\alpha\,n^2}{2+\alpha}^{(2+\alpha)n}
<
\bracket{8\,n}^{2(2+\alpha)n}.$$
Thus
$$\alpha>\frac{(1-\epsilon)\Delta-8}{4}.$$
Thus $\SL(G)\geq\quarter((1-\epsilon)\Delta-8)n$ for some $\Delta$-regular $n$-vertex graph $G$.
\end{proof}

Note that the lower bound in Theorem~\ref{thm:SlopeParameter} is within a factor of $2+\epsilon$ of the trivial upper bound $\SL(G)\leq\half\Delta n$. 

\section{Slope-Number: Proof of Theorem~\ref{thm:SlopeNumber}}
\label{sec:SlopeNumber}

In this section we extend the method developed in Section~\ref{sec:SlopeParameter} to prove a lower bound on the slope-number of graphs with bounded degree.

\begin{lemma}
\label{lem:SlopeNumber}
The number of labelled $n$-vertex $m$-edge graphs with slope-number at most $k$ is at most
$$\bfrac{50n^2(k+1)}{2n+k}^{2n+k}\binom{k(n-1)}{m}.$$
\end{lemma}

\begin{proof} Consider $V(G)=\{1,2,\dots,n\}$ for every labelled $n$-vertex $m$-edge graph $G$ with slope-number at most $k$. For every such graph $G$, fix a $k$-slope drawing of $G$ represented by a point set $\{(x_i(G),y_i(G)):1\leq i\leq n\}$ and slope set $\{s_\ell(G):1\leq\ell\leq k\}$. Thus for every edge $ij$ of $G$, the slope of the line through $(x_i(G),y_i(G))$ and $(x_j(G),y_j(G))$ equals $s_\ell(G)$ for some $\ell$. Without loss of generality, every $s_\ell(G)<\infty$. Define \PP\ as in the proof of Lemma~\ref{lem:Basic}. In addition, for all $i,j$ with $1\leq i<j\leq n$, define
$Q_{i,j}(G):=x_i(G)-x_j(G)$. Let $\QQ:=\{Q_{i,j}:1\leq i<j\leq n\}$. By Lemma~\ref{lem:MilnorThom} with $D=2$, $d=2n+k$, and $t=\binom{n}{2}(k+1)$, the number of sign patterns of $\PP\cup\QQ$ is at most $$\bfrac{50n^2(k+1)}{2n+k}^{2n+k}.$$ 

Fix a sign pattern $\sigma$ of $\PP\cup\QQ$. As in Lemma~\ref{lem:Basic}, from $\sigma$ restricted to \PP\ we can reconstruct the collinear subsets of vertices. Moreover, from $\sigma$ restricted to \QQ, we can reconstruct the order of the vertices within each collinear subset. Observe that at most $n-1$ edges have the same slope in a geometric drawing. Thus every $k$-slope graph representable by $\sigma$ is a subgraph of a fixed graph with at most $k(n-1)$ edges. Hence $\sigma$ corresponds to at most $\binom{k(n-1)}{m}$ labelled $k$-slope graphs on $m$ edges. The result follows.
\end{proof}

\begin{proof}[Proof of Theorem~\ref{thm:SlopeNumber}]
Suppose that for some $\Delta\geq5$ and for some integer $k$, every $\Delta$-regular graph has slope-number at most $k$. By Lemmas~\ref{lem:NumberRegular} and \ref{lem:SlopeNumber}, for all $n$, 
$$\bfrac{n}{3\Delta}^{\Delta n/2}
\leq
\bfrac{50n^2(k+1)}{2n+k}^{2n+k}\binom{k(n-1)}{\half\Delta n}.$$
Let $\epsilon=\epsilon(\Delta)>0$ be specified later. For large $n>n(k,\Delta,\epsilon)$ there is a constant $c=c(k,\Delta)$ such that 
$$\binom{k(n-1)}{\half\Delta n}\leq c^n= n^{n/\log_c n}<n^{\epsilon n}.$$ 
Thus
$$\bfrac{n}{3\Delta}^{\Delta n/2}
\leq
\bfrac{50n^2(k+1)}{2n+k}^{2n+k} n^{\epsilon n}
<
\bracket{25n(k+1)}^{(2+\epsilon)n+k}.$$
For large $n>n(\epsilon,k,\Delta)$ we have
$$
\bracket{25n(k+1)}^{(1-\epsilon)\Delta n/2}
<
\bfrac{n}{3\Delta}^{\Delta n/2}
<
\bracket{25n(k+1)}^{(2+\epsilon)n+k}.$$
Thus $(1-\epsilon)\Delta n<(4+\epsilon)n+2k$. Choose $\epsilon>0$ such that 
$(1-\epsilon)\Delta>4+\epsilon$. We obtain a contradiction for large $n>\frac{2k}{(1-\epsilon)\Delta-(4+\epsilon)}$. Thus there exists a $\Delta$-regular graph with slope-number greater than $k$
\end{proof}

It remains open whether slope-number is bounded by a constant for all degree-$3$ or degree-$4$ graphs.

\section*{Acknowledgement}

Thanks to Emo Welzl and the Theory of Combinatorial Algorithms group at ETH Z\"urich, and to Carsten Thomassen and the Technical University of Denmark for their generous hospitality which enabled this collaboration. Thanks to Vida Dujmovi\'c and Nick Wormald for fruitful discussions. 


\def\soft#1{\leavevmode\setbox0=\hbox{h}\dimen7=\ht0\advance \dimen7
  by-1ex\relax\if t#1\relax\rlap{\raise.6\dimen7
  \hbox{\kern.3ex\char'47}}#1\relax\else\if T#1\relax
  \rlap{\raise.5\dimen7\hbox{\kern1.3ex\char'47}}#1\relax \else\if
  d#1\relax\rlap{\raise.5\dimen7\hbox{\kern.9ex \char'47}}#1\relax\else\if
  D#1\relax\rlap{\raise.5\dimen7 \hbox{\kern1.4ex\char'47}}#1\relax\else\if
  l#1\relax \rlap{\raise.5\dimen7\hbox{\kern.4ex\char'47}}#1\relax \else\if
  L#1\relax\rlap{\raise.5\dimen7\hbox{\kern.7ex
  \char'47}}#1\relax\else\message{accent \string\soft \space #1 not
  defined!}#1\relax\fi\fi\fi\fi\fi\fi} \def\cprime{$'$}

\appendix
\section{Derivation of Lemma~\ref{lem:NumberRegular}}

Let $f(n,\Delta)$ denote the number of labelled $\Delta$-regular $n$-vertex graphs. The first asymptotic bounds on $f(n,\Delta)$ were independently determined by \citet{BC-JCTA78} and \citet{Wormald78}. Refining these results, \citet{McKay-AC85} proved that for all $\Delta$ with $1\leq\Delta<\frac{2}{9}n$, $$f(n,\Delta)=\frac{(\Delta n)!}{
	(\Delta n/2)! \;
	2^{\Delta n/2} \;	
	(\Delta!)^n \;
	\e^{(\Delta^2-1)/4-\Oh{\Delta^3/n}}}.$$
The version of Stirling's formula due to \citet{Robbins-AMM55} states that for all $t\geq1$, 
\begin{equation*}
\label{eqn:Stirling}
t!\,=\,\sqrt{2\pi t}\;\bfrac{t}{\e}^t\e^{r(t)},
\end{equation*}
where  $1/(12t+1)<r(t)<1/12t$. Thus, for some constant $c$,
\begin{align*}
f(n,\Delta)\;\geq\;
&\frac{\sqrt{2\pi\Delta n}\;
	\bfrac{\Delta n}{\e}^{\Delta n} \;
	\e^{r(\Delta n)} 
}{
	\sqrt{\pi\Delta n}\;
	\bfrac{\Delta n}{2\e}^{\Delta n/2}\;
	\e^{r(\Delta n/2)} \;
	2^{\Delta n/2}\;
	\Delta^{\Delta n}\;
	\e^{c\Delta^2}
}
\\\geq\;&
\sqrt{2}\bfrac{n}{\Delta}^{\Delta n/2}\;
/\exp\bracket{
	\frac{\Delta n}{2}
	-\frac{1}{12\Delta n}+\frac{1}{6\Delta n+1}
	+c\Delta^2}.
\end{align*}
With $n>200c\Delta$, we have 
$$\frac{\Delta n}{2}-\frac{1}{12\Delta n}+\frac{1}{6\Delta n+1}+c\Delta^2
<
\frac{102\,\Delta n}{200}.$$
Thus
\begin{align*}
f(n,\Delta)
>\bfrac{n}{\e^{1.02}\Delta}^{\Delta n/2}
>\bfrac{n}{3\Delta}^{\Delta n/2}
\end{align*}
Lemma~\ref{lem:NumberRegular} follows.

\section{Products of Binomials}

\begin{lemma}
\label{lem:Binomials}
Let $n$ and $t$ be positive integers. Let $x_1,x_2,\dots,x_t$ be nonnegative integers with each $x_i\leq n$. Let $a$ and $b$ be the unique integers such that
$\sum_ix_i=(t-b)a+b(a+1)$ and $0\leq b\leq t-1$. Then 
$$\prod_{i=1}^t\binom{n}{x_i}\leq\binom{n}{a}^{t-b}\binom{n}{a+1}^{b}.$$
\end{lemma}

\begin{proof}
Choose $x_1,x_2,\dots,x_t$ to maximise $\prod_i\binom{n}{x_i}$. Suppose on the contrary that two of the $x_i$ differ by at least two. Without loss of generality $x_1\geq x_2+2$. Let $x'_i:=x_i$ except for $x'_1:=x_1-1$ and  $x'_2:=x_2+1$. Thus $0\leq x'_i\leq n$. By assumption
$\prod_i\binom{n}{x_i}\geq\prod_i\binom{n}{x'_i}$. Hence $\binom{n}{x_1}\binom{n}{x_2}\geq\binom{n}{x_1-1}\binom{n}{x_2+1}$.
It follows that $x_1\leq x_2+1$. This contradiction proves that all pairs of the $x_i$ differ by at most one. Thus $\prod_i\binom{n}{x_i}$ is maximised when $t-b$ of the $x_i$ equal $a$, and $b$ of the $x_i$ equal $a+1$.
\end{proof}

\end{document}